\documentclass[11pt]{amsart}
\usepackage{graphicx,amssymb}

\def\baselinestretch{1}

\theoremstyle{plain}
\newtheorem{theorem}{Theorem}

\newtheorem{proposition}[theorem]{Proposition}

\def\TL{\qopname\relax n{\it TL}}

\begin{document}

\title[Dual presentation and linear basis of the Temperley-Lieb algebras]
{Dual presentation and linear basis\\
of the Temperley-Lieb algebras}

\author[E.-K.~Lee and S.~J.~Lee]
  {Eon-Kyung Lee and Sang Jin Lee}

\address
{Eon-Kyung Lee \\
Department of Applied Mathematics \\
Sejong University  \\
Seoul 143-747, Korea}
\email{\rm eonkyung@sejong.ac.kr}

\address
{Sang Jin Lee \\
Department of Mathematics \\
Konkuk University \\
Seoul 143-701, Korea}
\email{\rm sangjin@konkuk.ac.kr}



\begin{abstract}
The braid group $B_n$ maps homomorphically into
the Temperley-Lieb algebra $\TL_n$.
It was shown by Zinno that the homomorphic images of
simple elements arising from the dual presentation
of the braid group $B_n$
form a basis for the vector space underlying the
Temperley-Lieb algebra $\TL_n$.
In this paper, we establish that there is a dual presentation
of Temperley-Lieb algebras that corresponds to
the dual presentation of braid groups,
and then give a simple geometric proof for Zinno's theorem,
using the interpretation of simple elements
as non-crossing partitions.
\end{abstract}

\maketitle

\section{Introduction}
Since Jones~\cite{Jon83,Jon87}
discovered the Jones polynomial for links
by investigating representations of braid groups
into Hecke algebras and Temperley-Lieb algebras,
Temperley-Lieb algebras have played important roles
in the quantum invariants of links and 3-manifolds.
The Temperley-Lieb algebra $\TL_n$
is defined on non-invertible generators
$e_1,\ldots,e_{n-1}$ with the relations:
$e_ie_j=e_je_i$ for $|i-j|\ge 2$; $e_i^2=e_i$;
$e_ie_{i\pm 1}e_i=\tau e_i$
along with a complex number $\tau$.
It is well-known that the dimension of
$\TL_n$ is the $n$th Catalan number
$\mathcal C_n=\frac 1{n+1}{2n\choose n}$.
Setting $t$ such that $\tau^{-1}=2+t+t^{-1}$,
and then setting $h_i=(t+1)e_i-1$,
we get an alternative presentation
of $\TL_n$ with invertible generators $h_1,\ldots,h_{n-1}$
satisfying the relations:
\begin{eqnarray}
& & h_ih_j=h_jh_i\quad\mbox{if $|i-j|\ge 2$};\label{eqn:TLgenG1}\\
& & h_ih_{i+1}h_i=h_{i+1}h_ih_{i+1};\label{eqn:TLgenG2}\\
& & h_i^2=(t-1)h_i+t;\label{eqn:TLgenG3}\\
& & h_ih_{i+1}h_i+h_ih_{i+1}+h_{i+1}h_i+h_i+h_{i+1}+1=0.\label{eqn:TLgenG4}
\end{eqnarray}

The braid group $B_n$ is defined by the Artin presentation,
where the generators are $\sigma_1,\ldots,\sigma_{n-1}$
and the defining relations are
\begin{eqnarray*}
&&\sigma_i\sigma_j=\sigma_j\sigma_i
    \quad\text{if $|i-j|\ge 2$};\\
&&\sigma_i\sigma_{i+1}\sigma_i=\sigma_{i+1}\sigma_i\sigma_{i+1}
    \quad\text{for $i=1,\ldots,n-2$}.
\end{eqnarray*}
The braid group $B_n$ maps homomorphically
into the Temperley-Lieb algebra $\TL_n$
under $\pi:\sigma_i\mapsto h_i$.
There is another presentation~\cite{BKL98}
with generators $a_{ji}$ $(1\le i<j\le n)$
and defining relations
\begin{eqnarray*}
&&a_{lk}a_{ji}=a_{ji}a_{lk}\quad\mbox{if $(l-j)(l-i)(k-j)(k-i)>0$};\\
&&a_{kj}a_{ji}=a_{ji}a_{ki}=a_{ki}a_{kj}\quad\mbox{for $i<j<k$}.
\end{eqnarray*}
The generators $a_{ji}$'s are related to the $\sigma_i$'s
by
$$
a_{ji}=\sigma_{j-1}\sigma_{j-2}\cdots\sigma_{i+1}\sigma_i
\sigma_{i+1}^{-1}\cdots\sigma_{j-2}^{-1}\sigma_{j-1}^{-1}.
$$
Bessis~\cite{Bes01} showed that there is a similar presentation,
called the \emph{dual presentation},
for Artin groups of finite Coxeter type.

Both the Artin and dual presentations
of the braid group $B_n$ determine a \emph{Garside monoid,}
as defined by Dehornoy and Paris~\cite{DehPar99}, where
the \emph{simple elements} play important roles.
Nowadays, it becomes more and more popular to
describe simple elements arising from the dual presentation
via non-crossing partitions.
Non-crossing partitions are useful
in diverse areas~\cite{Bes01,Bra01,BC04,BDM02,Mcc06},
because they have beautiful combinatorial structures.

\begin{figure}
\includegraphics[scale=1.1]{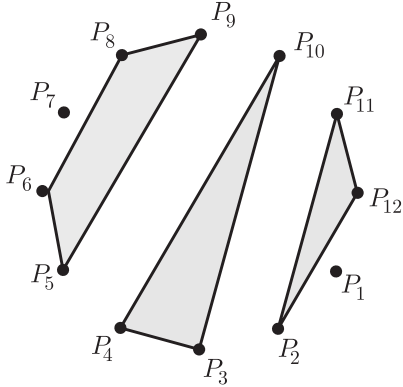}
\caption{The shaded regions show the blocks in the non-crossing
partition corresponding to the simple element
$[12,11,2]\,[10,4,3]\,[9,8,6,5]$ in $B_{12}$.}\label{fig:ncp}
\end{figure}

Let $P_1,\ldots,P_n$ be the points in
the complex plain given by $P_k=\exp(-\frac{2k\pi}ni)$.
See Figure~\ref{fig:ncp}.
Recall that a partition of a set is a collection of pairwise
disjoint subsets whose union is the entire set.
Those subsets (in the collection) are called blocks.
A partition of $\{P_1,\ldots,P_n\}$
is called a \emph{non-crossing partition}
if the convex hulls of the blocks are pairwise disjoint.

A positive word of the form
$a_{i_1i_2}a_{i_2i_3}\cdots a_{i_{k-1}i_k}$,
$i_1> i_2>\cdots> i_k$, is called a \emph{descending cycle}
and denoted $[i_1,i_2,\ldots,i_k]$.
Two descending cycles $[i_1,\ldots,i_k]$ and $[j_1,\ldots,j_l]$ are
said to be \emph{parallel} if the convex hulls
of $\{P_{i_1},\ldots,P_{i_k}\}$
and of $\{P_{j_1},\ldots,P_{j_l}\}$ are disjoint.
The simple elements are the products of parallel descending cycles.

We remark that the definition of simple elements
depends on the presentations.
For example, the simple elements arising from the Artin
presentation are in one-to-one correspondence with
permutations.
Throughout this note, we consider only
the simple elements arising from the dual presentation
of braid groups as above.

Note that simple elements are in one-to-one correspondence with
non-crossing partitions.
Our convention is that if a block in a non-crossing
partition consists of a single point, then
the corresponding descending cycle is the identity
(i.e.{} the descending cycle of length 0).
In particular, the number of the simple elements is
the $n$th Catalan number $\mathcal C_n$,
which is the dimension of $\TL_n$.
Zinno~\cite{Zin02} established the following result.

\begin{theorem}[Zinno's theorem]\label{thm:TL-Generator}
The homomorphic images of the simple elements arising
from the dual presentation of $B_n$
form a linear basis for the Temperley-Lieb algebra
$\TL_n$.
\end{theorem}

We explain briefly Zinno's proof.
It is known that the \emph{ordered reduced words}
$$
(h_{j_1}h_{j_1-1}\cdots h_{k_1})
(h_{j_2}h_{j_2-1}\cdots h_{k_2})\cdots
(h_{j_p}h_{j_p-1}\cdots h_{k_p}),
$$
where $j_i\ge k_i$, $j_{i+1}>j_i$ and $k_{i+1}>k_i$,
form a linear basis of $\TL_n$,
and Zinno showed that the matrix for writing the images
of simple elements as the linear combination of
the ordered reduced words is invertible.
Because the number of the simple elements is equal to
the dimension of $\TL_n$,
this proves the theorem.

In this note, we first establish
that there is a dual presentation of $\TL_n$.
We are grateful to David Bessis for pointing out that
the relation~(\ref{eqn:TLgenG4}) in the Temperley-Lieb algebra
presentation is equivalent to the forth relation
in the dual presentation in the following theorem.

\begin{theorem}[dual presentation of $\TL_n$]\label{thm:TL-DualPreA}
The Temperley-Lieb algebra $\TL_n$ has a presentation
with invertible generators
$g_{ji}$ $(1\le i<j\le n)$
satisfying the relations:
$$
\begin{array}{l}
g_{lk}g_{ji}=g_{ji}g_{lk}
    \quad\mbox{if\/ $(l-j)(l-i)(k-j)(k-i)>0$};\\
g_{kj}g_{ji}=g_{ji}g_{ki}=g_{ki}g_{kj}
    \quad\mbox{for $i<j<k$};\\
g_{ji}^2=(t-1)g_{ji}+t
    \quad\mbox{for $i<j$};\\
g_{ji}g_{kj}+tg_{kj}g_{ji}+g_{kj}+g_{ji}+tg_{ki}+1=0
    \quad\mbox{for $i<j<k$}.
\end{array}
$$
The new generators are related to the old ones
by
$$g_{ji}=h_{j-1}h_{j-2}\cdots h_{i+1}h_i
h_{i+1}^{-1}\cdots h_{j-2}^{-1}h_{j-1}^{-1}.$$
\end{theorem}

Using the above presentation, we give a new proof of Zinno's theorem
in \S\ref{sec:dimension}.
We exploit non-crossing partitions so as to make
the proof easy and intuitive.
For the proof, we show that
any monomial in the $h_i^{\pm 1}$'s can be written as
a linear combination of the images of simple elements.
Therefore the images of simple elements span $\TL_n$.
As a result, they form a linear basis of $\TL_n$
because the number of simple elements is equal to
the dimension of $\TL_n$.

We remark that it seems possible to prove
the linear independence of
the images of the simple elements
directly from the relations in the dual
presentation of $\TL_n$
(without using the fact that the dimension
of $\TL_n$ is the same as the number of simple elements),
but that would be beyond the scope of this note
because it would require repeating all the arguments
used in the proof for the embedding of the positive braid monoid
in the braid group.

\subsection*{Acknowledgements}
We are very grateful to David Bessis
for the intensive discussions
during his visit to Korea Institute for Advanced
Study in June 2003.

\section{Dual presentation of the Temperley-Lieb algebras}
\label{sec:DualPre}

Let $D^2$ be the disc in the complex plane with radius 2
and $P_1,\ldots,P_n$ be the points in $D^2$
given by $P_k=\exp(-\frac{2k\pi}ni)$.
Let $D_n=D^2\setminus\{P_1,\ldots,P_n\}$.
The braid group $B_n$ can be regarded
as the group of self-homeomorphisms of $D_n$ that
fix the boundary pointwise, modulo isotopy relative to the
boundary. The generators $\sigma_i$ and $a_{ji}$ correspond
to the positive half Dehn-twists along the arcs
${P_iP_{i+1}}$ and ${P_iP_j}$, respectively.

\begin{figure}
\begin{tabular}{ccc}
\includegraphics[scale=1]{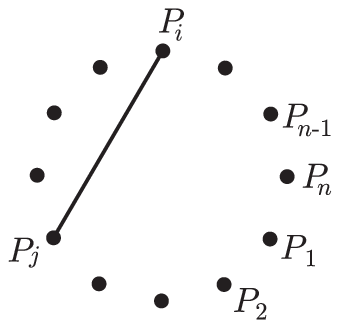}\quad\mbox{}&\mbox{}\quad
\includegraphics[scale=1]{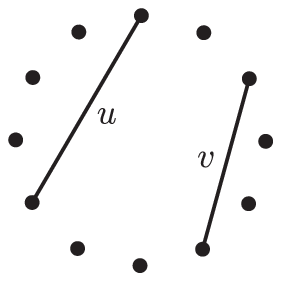}\quad\mbox{}&\mbox{}\quad
\includegraphics[scale=1]{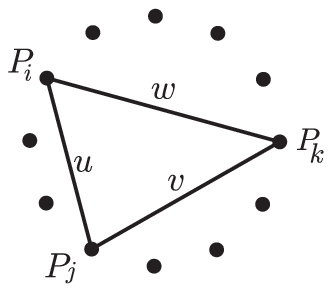}\\
(a) & (b) & (c)
\end{tabular}
\caption{}\label{fig:DualGen}
\end{figure}

Let $\mathcal L_n=\{P_iP_j\mid i\ne j\}$
be the set of line segments as in Figure~\ref{fig:DualGen}~(a).
We say that a pair $(u,v)\in\mathcal L_n^2$ is \emph{parallel}
if $u$ and $v$ are disjoint as in Figure~\ref{fig:DualGen}~(b),
and \emph{admissible} if $u=P_iP_j$ and $v=P_jP_k$ for some
pairwise distinct points $P_i$, $P_j$ and $P_k$
which are in counterclockwise order on the unit circle
as in Figure~\ref{fig:DualGen}~(c).
A triple $(u,v,w)\in\mathcal L_n^3$ is said to be
\emph{admissible} if so are all the pairs
$(u,v)$, $(v,w)$ and $(w,u)$.
The dual presentation of $B_n$ can be written as follows:
$$
B_n=\biggl\langle
{a_u\atop (u\in\mathcal L_n)}
\Bigm|
\begin{array}{l}
a_ua_v=a_va_u\quad\mbox{if $u$ and $v$ are parallel}\\
a_ua_v=a_va_w=a_wa_u\quad\mbox{if $(u,v,w)$ is admissible}
\end{array}
\bigg\rangle.
$$
It is easy to see the following:
(i) if $u\cap v=\{P_i\}$ for some $P_i$ in the unit circle
then exactly one of $(u,v)$ and $(v,u)$ is admissible;
(ii) if $(u,v)$ is admissible, then
$a_ua_v$ can be written in three ways as in the presentation,
but $a_va_u$ is not equivalent to any other positive word
on the $a_u$'s;
(iii) $a_ua_v$ is a simple element if and only if $(u,v)$ is parallel
or admissible.

Now we prove Theorem~\ref{thm:TL-DualPreA}.
The theorem can be rewritten as follows.
Its proof is elementary.
However, we present it for completeness.

\begin{theorem}[dual presentation of $\TL_n$]\label{thm:TL-DualPre}
The Temperley-Lieb algebra $\TL_n$ has a presentation
with invertible generators
$g_u$ $(u\in \mathcal L_n)$
satisfying the relations:
\begin{eqnarray}
&&g_ug_v=g_vg_u\quad\mbox{if\/ $u$ and $v$ are parallel};
    \label{eqn:TLgenD1}\\
&&g_ug_v=g_vg_w=g_wg_u\quad\mbox{if\/ $(u,v,w)$ is admissible};
    \label{eqn:TLgenD2}\\
&&g_u^2=(t-1)g_u+t\quad\mbox{for $u\in\mathcal L_n$};
    \label{eqn:TLgenD3}\\
&&g_vg_u+tg_ug_v+g_u+g_v+tg_w+1=0\quad\mbox{if\/ $(u, v,w)$ is admissible}.
    \label{eqn:TLgenD4}
\end{eqnarray}
\end{theorem}

\begin{proof}
From the results on the dual presentation of $B_n$ in ~\cite{BKL98},
it follows that the relations~(\ref{eqn:TLgenG1}) and (\ref{eqn:TLgenG2})
are equivalent to the relations~(\ref{eqn:TLgenD1}) and (\ref{eqn:TLgenD2}).

\medskip
Assume the relations (\ref{eqn:TLgenG1}), (\ref{eqn:TLgenG2}),
and hence (\ref{eqn:TLgenD1}), (\ref{eqn:TLgenD2}).

(\ref{eqn:TLgenD3}) $\Rightarrow$ (\ref{eqn:TLgenG3}) \ \
It is clear since (\ref{eqn:TLgenG3}) is a special case of (\ref{eqn:TLgenD3}).

(\ref{eqn:TLgenG3}) $\Rightarrow$ (\ref{eqn:TLgenD3}) \ \
It is clear since each $g_u$ is conjugate to $h_i$ (for some $i$)
by a monomial in the $h_j$'s.

\medskip
Now assume the relations (\ref{eqn:TLgenG1}), (\ref{eqn:TLgenG2}),
(\ref{eqn:TLgenG3}), and hence (\ref{eqn:TLgenD1}), (\ref{eqn:TLgenD2}),
(\ref{eqn:TLgenD3})

(\ref{eqn:TLgenD4}) $\Rightarrow$ (\ref{eqn:TLgenG4}) \ \
Let $u={P_{i+2}P_{i+1}}$, $v={P_{i+1}P_i}$ and $w={P_iP_{i+2}}$.
Then $g_{u}=h_{i+1}$, $g_{v}=h_i$ and $(u,v,w)$ is admissible.
Since
$h_ih_{i+1}h_i=g_vg_ug_v
=g_{v}g_{v}g_{w}=((t-1)g_{v}+t)g_{w}
=(t-1)g_{v}g_{w}+tg_{w}=(t-1)g_{u}g_{v} +tg_{w}$,
\begin{eqnarray*}
\lefteqn{h_ih_{i+1}h_i+h_ih_{i+1}+h_{i+1}h_i+h_i+h_{i+1}+1}\\
&=&\left((t-1)g_{u}g_{v} +tg_{w}\right) + g_{v}g_{u}+g_{u}g_{v}+g_{v}+g_{u}+1\\
&=& g_{v}g_{u}+tg_{u}g_{v}+g_{u}+g_{v}+tg_{w}+1 \ = \ 0.
\end{eqnarray*}

(\ref{eqn:TLgenG4}) $\Rightarrow$ (\ref{eqn:TLgenD4}) \ \
Note that for each admissible triple $(u',v',w')$,
there is a self-homeomorphism of $D_n$ sending $(u',v',w')$
to $(u,v,w)$. Therefore, there is a monomial $x$
in the $h_j$'s such that
$x g_{u'}x^{-1}=g_{u}$, $x g_{v'}x^{-1}=g_{v}$ and
$x g_{w'}x^{-1}=g_{w}$, simultaneously.
Let $u'={P_{i+2}P_{i+1}}$, $v'={P_{i+1}P_i}$ and $w'={P_iP_{i+2}}$.
In the same way as in
(\ref{eqn:TLgenD4}) $\Rightarrow$ (\ref{eqn:TLgenG4}),
we obtain
\begin{eqnarray*}
\lefteqn{g_{v}g_{u}+tg_{u}g_{v}+g_{u}+g_{v}+tg_{w}+1} \\
&=& x \left( g_{v'}g_{u'}+tg_{u'}g_{v'}+g_{u'}+g_{v'}
    +tg_{w'}+1 \right) x^{-1} \\
&=& x \left( h_ih_{i+1}h_i
   +h_ih_{i+1}+h_{i+1}h_i
   +h_i+h_{i+1}+1 \right) x^{-1} \\
&=& 0.
\end{eqnarray*}
\vskip-\baselinestretch\baselineskip
\end{proof}

\section{A new proof of Zinno's theorem}
\label{sec:dimension}

Before starting the proof of Zinno's theorem,
let us observe the relations
$g_u^2=(t-1)g_u+t$ and $g_vg_u+tg_ug_v+g_u+g_v+tg_w+1=0$
in the dual presentation of $\TL_n$.
Among the monomials in the relations,
all except $g_u^2$ and $g_vg_u$ are images of simple elements.
Therefore the relations can be interpreted
as instructions for converting a product
of two generators into a linear combination of
the images of simple elements:
\begin{eqnarray*}
g_u^2&=&(t-1)\pi(a_u)+t;\\
g_vg_u&=&-t\pi(a_ua_v) -\pi(a_u) -\pi(a_v) -t\pi(a_w)-1.
\end{eqnarray*}

Generalizing this idea,
we will show in Proposition~\ref{thm:span}
that for a simple element $A$ and an Artin generator $\sigma_i$,
the homomorphic image $\pi(A\sigma_i)$ in $\TL_n$
can be written as a linear combination
of the images of simple elements.

Recall that the simple elements are in one-to-one correspondence
with non-crossing partitions.
For a simple element $A$,
take union of the convex hulls of the blocks
in the non-crossing partition of $A$,
and then remove those containing only one point.
The resulting set is called the \emph{underlying space} of $A$
and denoted $\bar A$.

\begin{figure}
\begin{tabular}{cc}
\includegraphics[scale=1.1]{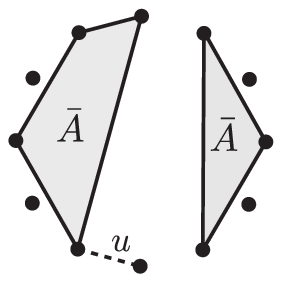}
\ \ \mbox{}\raisebox{14mm}{\Large $=$}\ \
\includegraphics[scale=1.1]{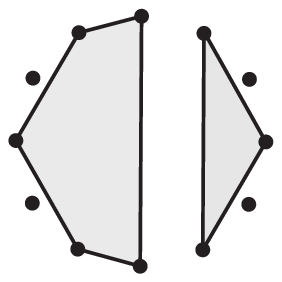}\\[5mm]
\includegraphics[scale=1.1]{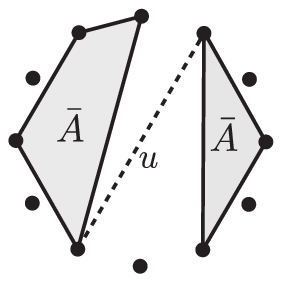}
\ \ \mbox{}\raisebox{14mm}{\Large $=$}\ \
\includegraphics[scale=1.1]{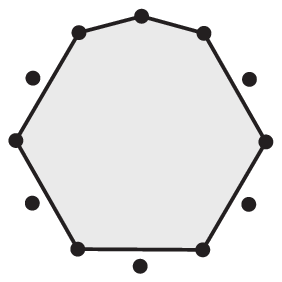}
\end{tabular}
\caption{In the left hand sides $\bar A$ and $u$ are depicted
as shaded regions and dotted lines.
The right hand sides show the underlying spaces of $Aa_u$'s.}
\label{fig:Simple}
\end{figure}

\begin{figure}
$$
\includegraphics[scale=1.1]{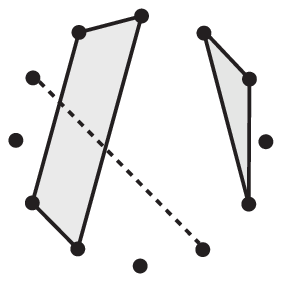}\qquad
\includegraphics[scale=1.1]{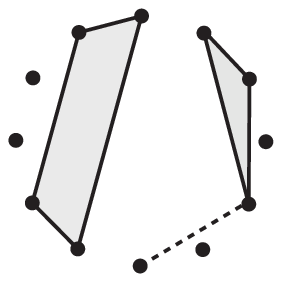}\qquad
\includegraphics[scale=1.1]{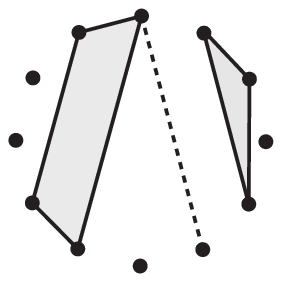}
$$
\caption{The shaded regions and the dotted lines represent
the underling space of a simple element $A$ and
an element $u\in\mathcal L_n$, respectively.
In this case, $A a_u$ is not a simple element.}
\label{fig:SimpleNot}
\end{figure}

It is known that for a simple element $A$ and $u\in\mathcal L_n$,
$Aa_u$ is a simple element if and only if
for any $w\in \mathcal L_n$ with $w\subset \bar A$,
the product $a_wa_u$ is a simple element, in other words,
$(w,u)$ is parallel or admissible~\cite[Corollary 3.6]{BKL98}.
Figure~\ref{fig:Simple} shows typical cases of $(\bar A, u)$
such that $Aa_u$ becomes a simple element,
and Figure~\ref{fig:SimpleNot} shows some cases of $(\bar A,u)$
such that $Aa_u$ is not a simple element.

It is easy to see that if $\bar A$ and $u$ satisfy one
of the following conditions,
then $Aa_u$ is a simple element and its underlying
space is the union of the convex hulls of components of $\bar A\cup u$.

\begin{itemize}
\item $\bar A$ and $u$ are disjoint.
\item $\bar A$ and $u$ intersect at the boundary of $u$
as in the left hand sides of Figure~\ref{fig:Simple}.
Intuitively, when we stand at an intersection point,
with $u$ on the right and the component of $\bar A$ containing
the intersection point on the left,
we become to face towards the inside of the unit circle.
\end{itemize}

\begin{proposition}\label{thm:span}
For a simple element $A$ and an Artin generator $\sigma_i$,
$\pi(A\sigma_i)$ can be expressed as a linear combination
of the images of simple elements.
\end{proposition}

\begin{proof}
Let $u={P_iP_{i+1}}$.
Then $\sigma_i=a_u$ and $\pi(\sigma_i)=g_u$.
We prove the assertion in three cases.

\smallskip
\noindent\textsl{Case 1.\ }
If $u\subset\bar A$, then $\bar A$ and $u$ are
as in Figure~\ref{fig:case1}~(a).
Let $B$ be the simple element whose underlying space
is as in Figure~\ref{fig:case1}~(b).
More precisely, the non-crossing partition of $B$ is obtained
from that of $A$ by making $\{P_i\}$ a new block.
Then $A=Ba_u$ and
\begin{eqnarray*}
\pi(Aa_u)
&=&\pi(Ba_u^2)=\pi(B)g_u^2=\pi(B)((t-1)g_u+t)\\
&=&(t-1)\pi(Ba_u)+t\pi(B)=(t-1)\pi(A)+t\pi(B).
\end{eqnarray*}

\begin{figure}\tabcolsep=30pt
\begin{tabular}{cc}
\includegraphics[scale=1]{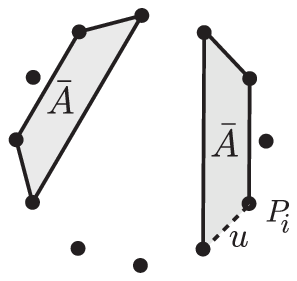}&
\includegraphics[scale=1]{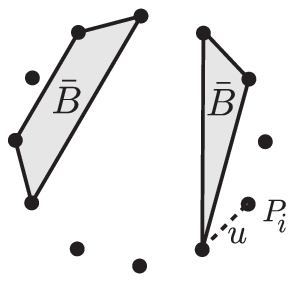}\\
(a) & (b)
\end{tabular}
\caption{$A=Ba_u$ if $\bar A$, $u$ and $\bar B$ are as above.}
\label{fig:case1}
\end{figure}

\smallskip
\noindent\textsl{Case 2.\ }
If $u\not\subset\bar A$ and $P_i\not\in\bar A$,
then $\bar A$ and $u$ are as in Figure~\ref{fig:case2}.
In this case, $Aa_u$ itself is a simple element.
\begin{figure}\tabcolsep=30pt
\begin{tabular}{cc}
\includegraphics[scale=1]{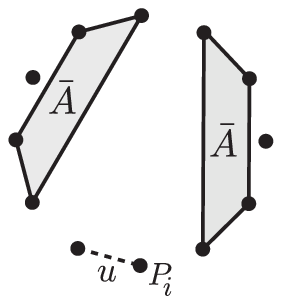}&
\includegraphics[scale=1]{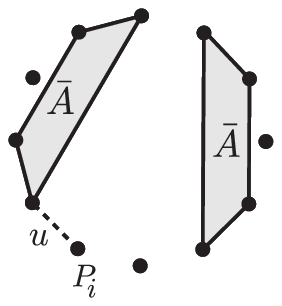}
\end{tabular}
\caption{$Aa_u$ is a simple element.}\label{fig:case2}
\end{figure}

\smallskip
\noindent\textsl{Case 3.\ }
If $u\not\subset\bar A$ and $P_i\in \bar A$, then
$\bar A$ and $u$ are either (a) or (b)
of Figure~\ref{fig:case3}, depending on whether $P_{i+1}$ belongs to
$\bar A$ or not.
Let $v$ be the line segment containing $P_i$
such that $Ba_v=A$ for some simple element $B$ as in (c) and (d)
of Figure~\ref{fig:case3}.
(More precisely, $v=P_iP_j$ for some $P_j$
such that $P_iP_j\subset \bar A$ and
the interior of $P_{i+1}P_j$ does not intersect $\bar A$.)
Let $w$ be the line segment connecting the endpoints of $u$
and $v$ other than $P_i$.
\begin{figure}\tabcolsep=30pt
\begin{tabular}{cccc}
\includegraphics[scale=1]{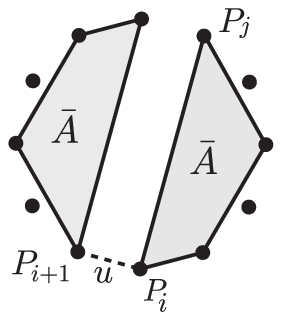}&
\includegraphics[scale=1]{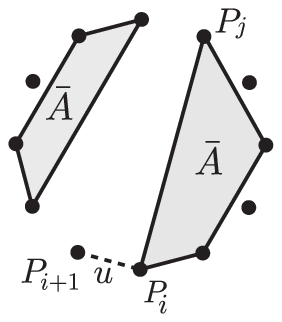}\\
(a) & (b)\\[5mm]
\includegraphics[scale=1]{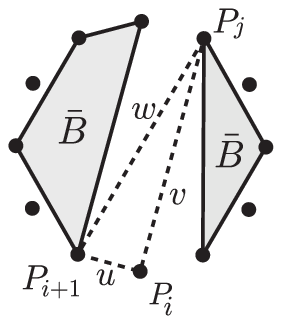}&
\includegraphics[scale=1]{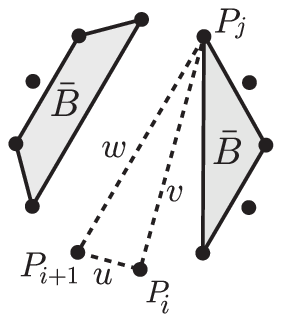}\\
(c) & (d)
\end{tabular}
\caption{}\label{fig:case3}
\end{figure}
Then $(u,v,w)$ is admissible and
\begin{eqnarray*}
\pi(Aa_u)
&=&\pi(Ba_va_u)=\pi(B)g_vg_u\\
&=&-\pi(B)(tg_ug_v+g_u+g_v+tg_w+1)\\
&=&-t\pi(Ba_ua_v)-\pi(Ba_u)-\pi(Ba_v)-t\pi(Ba_w)-\pi(B).
\end{eqnarray*}
Note that $Ba_ua_v$, $Ba_u$, $Ba_v$ and $Ba_w$ are simple elements.
\end{proof}

\def\temp{
\begin{proof}[Proof of Theorem \ref{thm:TL-Generator}]
Let $V_n$ be the subspace (of $\TL_n$) spanned by the images of simple elements.
Using Proposition~\ref{thm:span} and induction on the word length
of monomials in the $\sigma_i$'s, it is easy to show that
the images of monomials in the $\sigma_i$'s belong to $V_n$.
This implies that the images of $n$-braids belong to $V_n$,
because $h_i^{-1}=\frac1t h_i-\frac{t-1}t$.
Therefore $V_n=\TL_n$.
Because the number of simple elements is equal to the dimension of $\TL_n$,
this proves that the images of simple elements form
a linear basis of $\TL_n$
\end{proof}
}

\begin{proof}[Proof of Theorem \ref{thm:TL-Generator}]
Let $V_n$ be the subspace (of $\TL_n$) spanned by the images of simple elements.
Since the number of simple elements is equal to the dimension of $\TL_n$,
the images of simple elements form a linear basis of $\TL_n$
if we show that $V_n = \TL_n$ (i.e.
every monomial in the $h_i^{\pm 1}$'s belongs to $V_n$).

Observe that $h_i^{-1}=\frac1t h_i-\frac{t-1}t$ and $h_i = \pi(\sigma_i)$
for all $i$.
Therefore, it suffices to show that
the images of monomials in the $\sigma_i$'s belong to $V_n$.
Use induction on the word length of monomials in the $\sigma_i$'s.
By Proposition~\ref{thm:span}, it is easy to get the desired result.
\end{proof}

\end{document}